\documentclass{article}
\usepackage{amsmath}
\usepackage{amsthm}
\usepackage{amssymb}
\usepackage{graphicx}
\usepackage{epsfig}
\def\ml{l\kern-0.035cm\char39\kern-0.03cm}

\theoremstyle{plain}
\newtheorem{theorem}{Theorem}
\newtheorem{Theorem}[theorem]{Theorem}
\newtheorem{Lemma}[theorem]{Lemma}
\newtheorem{Corollary}[theorem]{Corollary}

\newtheorem{Observation}[theorem]{Observation}

\theoremstyle{definition}

\newtheorem{Question}[theorem]{Question}

\theoremstyle{remark}

\begin{document}

\title{\bf The Thue choice number versus the Thue chromatic number of graphs}

\author{Erika \v Skrabu\ml\'akov\'a \\
\small Faculty BERG\\[-0.8ex]
\small Technical University of Ko\v sice\\[-0.8ex] 
\small Ko\v sice, Slovakia\\
\small\tt erika.skrabulakova@tuke.sk}


\maketitle

\begin{abstract}
We say that a vertex colouring $\varphi$ of a graph $G$ is nonrepetitive if there is no positive integer $n$ and a path on $2n$ vertices $v_{1}\ldots v_{2n}$ in $G$ such that the associated sequence of colours $\varphi(v_{1})\ldots\varphi(v_{2n})$ satisfy $\varphi(v_{i})=\varphi(v_{i+n})$ for all $i=1,2,\dots,n$. 
The minimum number of colours in a nonrepetitive vertex colouring of $G$ is  the \emph{Thue chromatic number} $\pi (G)$. For the case of vertex list colourings  the {\em Thue choice number} $\pi_{\l}(G)$ of $G$ denotes the smallest integer $k$ such that for every list assignment 
$L:V(G)\rightarrow 2^{\mathbb{N}}$ with minimum list length at least $k$, there is a nonrepetitive vertex colouring of $G$ from the assigned lists. 
Recently it was proved that the Thue chromatic number and the Thue choice number of the same graph may have an arbitrary large difference in some classes of graphs. 
Here we give an overview of the known results where we compare these two parameters for several families of graphs and we also give a list of open problems on this topic.

  

\end{abstract}

\section{Introduction}
\label{intro}
A sequence $S$ is called {\em repetitive} if it contains a subsequence of consecutive terms $r_1 r_2 \dots r_{2n}$ called a {\em repetition}, e. g. for which it holds $r_i=r_{i+n}$ for all $i\in\{1,2,\dots,n\}$. Else $S$ is called  {\em nonrepetitive} (or  {\em square-free}). The words ABBA, BARBAR, HOTSHOTS are examples of repetitive sequences while MINIMIZE, NUMBER or COLOURS represents the nonrepetitive ones. \\
Nonrepetitive sequences were first studied by Norwegian mathematician Axel Thue at the beginning of the last century. Via his investigation of word structures (see \cite{Be95}, \cite{Be92}, \cite{Th06}) he became one of the ``fathers'' of a new branch of combinatorics - {\em Combinatorics on Words}, see \cite{Lo83}. In his famous paper from 1906 \cite{Th06} he showed the existence of arbitrarily long nonrepetitive sequences over three digit alphabet. 
Nonrepetitive sequences found their applications in many different areas but only in 80's they appeared also in graph theory \cite{Cur87} and with connections to graph colourings only in 2002 via the seminal paper of Alon et. al. \cite{AGHR02}. \\ 
A nonrepetitive colouring of graph is defined as follows: Let $\varphi$ be a colouring of the vertices of a graph $G$.  We say that $\varphi$ is a {\em nonrepetitive vertex colouring} of $G$ if for every path on $2n$ vertices $v_1$, $v_2, \dots, v_{2n}$ in $G$ the associated sequence of colours $\varphi(v_1)$ $\varphi(v_2) \dots \varphi(v_{2n})$ is not a repetition. The minimum number of colours in a nonrepetitive vertex colouring of a graph $G$ is the {\em Thue chromatic number} $\pi(G)$. For the case of list colourings the {\em Thue choice number} $\pi_{\l}(G)$ of a graph $G$ denotes the smallest integer $k$ such that for every list assignment $L: V(G)\to 2\sp{\mathbb{N}}$ with minimum list length at least $k$, there is a colouring of the vertices from the assigned lists such that the sequence of vertex colours of no path in $G$ forms a repetition. \\
It is known that in general the chromatic number and the choice number of the same graph may have an arbitrary large difference - consider for instance the class of complete bipartite graphs where there is no constant bound on the choice number (see \cite{ERT79}).
A similar result was for longer time not known for the case of the Thue chromatic number and the Thue choice number. Recently Fiorenzi et al. \cite{FOOMZ11} proved that the Thue chromatic number and the Thue choice number of the same graph may have an arbitrary large difference in some classes of graphs. On the other hand, there exist families of graphs where one can write the symbol of equality between these two parameters for every graph from the family. As both of the parameters are widely studied in the last years, the number of results grows very fast. Therefore, the purpose of this paper is to survey the results and problems on this topic. 

\section{General bounds}
\label{sec:1}

\subsection{Notations and basic observations} \label{Not}

The name {\em Thue number}\footnote{in the present called the Thue chromatic index and abbreviated $\pi'(G)$} and abbreviation $\pi(G)$ for the graph parameter dealing with nonrepetitive sequences for the first time appeared in the paper of Alon et al. \cite{AGHR02}. Unfortunatelly, in connection with edge colourings. The vertex version parameter\footnote{nowdays called the Thue chromatic number and abbreviated $\pi(G)$} was called here the {\em vertex Thue number} with no abbreviation at all. This was the cause of a lot of misunderstanding, as many papers in this area employ various notations and terminologies. Although some of the authors still stick to the original terminology of Alon et al. \cite{AGHR02}, many others, in order to get better transparence to the terminology, agreed with the following notations: each  Thue graph parameter  connected with nonrepetitive edge colourings will be called {\em index} and abbreviated with single quotation mark - apostrophe 
(e.g. the Thue chromatic index, $\pi'(G)$ - see \cite{BGKNP07}, \cite{Gr07}; 
the Thue choice index, $\pi_{\l}'(G)$ - see \cite{GoMoPi14}; 
the facial Thue chromatic index, $\pi'_f(G)$ - see \cite{GoMoPi14}, \cite{HJSS09}, \cite{JeSk09}, \cite{ES09};
the facial Thue choice index, $\pi_{f\l}'(G)$ - see \cite{GoMoPi14}, \cite{JP12}, \cite{JSES11}), 
while the Thue graph parameter connected with nonrepetitive vertex  colourings or total colourings will be called {\em number} and abbreviated without apostrophe
(e.g. the Thue chromatic number, $\pi(G)$ - see \cite{BGKNP07}, \cite{CLSMMS11}, \cite{CDS12}, \cite{GoMoPi14}, \cite{Gr06}, \cite{Gr07a}, \cite{Gr07b}, \cite{Gr07}, \cite{HaJe12}, \cite{HSS11}, \cite{KPZ13}, \cite{KSX12}, \cite{PSST14}, \cite{PeZm09},  \cite{SchUm02}; 
the Thue choice number, $\pi_{\l}(G)$ - see \cite{FOOMZ11}, \cite{GoMoPi14}, \cite{GPZ10}\footnote{here abbreviated $\pi_{ch}(G)$}, \cite{KM13};   
the facial Thue choice number, $\pi_{f\l}(G)$ - see \cite{GoMoPi14}, \cite{PSS14}; 
the facial Thue chromatic number, $\pi_f(G)$ - see \cite{BaCz11}, \cite{GoMoPi14}, \cite{HaJe11}, \cite{HaJe12},  unhapilly, with the same abbreviation like 
the fractional Thue chromatic number, $\pi_f(G)$ -  see \cite{KPZ13},  \cite{ZZ13}; for the Thue parameters related to total Thue colourings see \cite{KLMS15}, \cite{JSES13}). \\
We will also follow this idea and use the abbreviation $\pi(G)$ for the Thue chromatic number and $\pi_{\l}(G)$ for the Thue choice number of a graph $G$. Except of the few notation defined throughout the paper we will use the standard terminology according to Bondy and Murty \cite{BoMa08}. The maximum degree of a graph $G=G(V,E)$ will be denoted by $\Delta$ and order of the graph $|V(G)|=n$. \\

A vertex colouring of a graph $G$ such that no two adjacent vertices receive the same colour is called a {\em proper vertex coloring}.  The minimum number of colours in a proper vertex colouring of a graph $G$ is the {\em chromatic number} of $G$, $\chi(G)$. As adjacent vertices receive distinct colours in every nonrepetitive colouring, it is trivially proper. A proper colouring with no $2$-coloured $P_4$ is called a {\em star colouring} since each bichromatic subgraph is a star forest. The {\em star chromatic number}, $\chi_{st}(G)$, is the minimum number of colours in a star colouring of $G$ (see \cite{ACKKR04}). The relation between the chromatic number of a graph $G$, its star chromatic number and Thue chromatic number can be expressed as follows: \\

\begin{Observation}
$\chi(G)\leq\chi_{st}(G)\leq\pi(G)$. \\
\end{Observation} 

The basic observation on nonrepetitive vertex colouring is the following: As every nonrepetitive $k$-colouring of $G$ can be considered as a nonrepetitive list colouring of $G$ from identical lists of size $k$, the Thue  choice number of $G$ is a natural upper bound for the Thue  chromatic number of $G$. \\

\begin{Observation} 
$\pi(G)\leq \pi_{\l}(G)$. \\
\end{Observation} 

Another simple observation is that the bounds on the Thue chromatic number achieved by a probabilistic approach also hold for the Thue choice number.

\subsection{General bounds for $\pi(G)$ and $\pi_{\l}(G)$ based on $\Delta$} 

Since 2002 it is known that graphs with maximum degree $\Delta$ are nonrepetitively $O(\Delta\sp{2})$-colourable. The first upper bound for $\pi(G)$ in the form $c \Delta^2$ comes from the remark on vertex colourings in the seminal paper of Alon et al. \cite{AGHR02} and it can be achieved by probabilistic method using Lov\'asz local lemma in the proof (see also e. g.  \cite{AlSp00}, \cite{ELLL75}, \cite{HSS11}, \cite{KSX12}, \cite{MoLLL09}, \cite{Peg09}, \cite{JSES11}). \\

\begin{Theorem} \label{ckratdelta} 
(Alon, Grytczuk, Ha\l uszczak, Riordan, 2002, \cite{AGHR02}) \\
There exists a constant $c> 0$ such that
$\pi(G)\leq c \Delta\sp{2}$, for all graphs $G$ with maximum degree $\Delta$. \\
\end{Theorem} 

Alon et al. \cite{AGHR02} were dealing also with the lower bound for the parameter $\pi(G)$ and showed the following: \\

\begin{Theorem} \label{lower} 
(Alon, Grytczuk, Ha\l uszczak, Riordan, 2002, \cite{AGHR02}) \\
There exists a constant $c> 0$ with the following property: For every integer $\Delta > 1$, there exists a graph $G$ with maximum degree $\Delta$ such that every nonrepetitive vertex colouring of $G$ uses at least $c\frac{\Delta\sp{2}}{\log \Delta}$ colours. \\
\end{Theorem}

The proof is probabilistic, hence, it is also valid for $\pi_{\l}(G)$ and we can sum up Theorem~\ref{ckratdelta} and Theorem~\ref{lower} as follows: \\

\begin{Theorem} 
(Alon, Grytczuk, Ha\l uszczak, Riordan, 2002, \cite{AGHR02}) \\
$c_1\frac{\Delta\sp{2}}{\log \Delta}\leq \pi(G) \leq \pi_{\l}(G)\leq c_2\cdot \Delta\sp{2}$ for some constants $c_1$ and $c_2$. \\
\end{Theorem} 

The originally proved constant $c_2=2e\sp{16}$, was improved by Grytczuk to 36 \cite{Gr06} and 16 \cite{Gr07a}, later by Kolipaka et al. \cite{KSX12} to 10,4. These proofs are based on Lov\'asz local lemma, hence valid for $\pi_{\l}(G)$ too. \\
Harant and Jendro\ml \ \cite{HaJe12} for graphs with maximum degree $\Delta \geq 2$ proved that $\pi(G)\leq \pi_l(G) \leq \lceil 12,92(\Delta-1)\sp{2}\rceil$. 

Dujmovi\'c et al. \cite{DJKW13} using the entropy compression method (see e.g.  \cite{GoMoPi14}, \cite{GKM10}, \cite{MoLLL09},  \cite{MoTa10},\cite{JP12},  \cite{PSS14}) also improved the constant $c$ in the upper bound $c \Delta\sp{2}$ and showed that for large graphs $c$ even tends to 1:  \\

\begin{Theorem} \label{dujmovic} 
(Dujmovi\'c, Joret, Kozik, Wood, 2015+, \cite{DJKW13}) \\
For every graph $G$ with maximum degree $\Delta>1$, 
$$\pi(G) \leq \pi_{\l}(G)\leq\left\lceil \left( 1+\frac{1}{\Delta\sp{\frac{1}{3}}-1}+\frac{1}{\Delta\sp{\frac{1}{3}}} \right)\Delta\sp{2}\right\rceil=\Delta\sp{2}+2\Delta\sp{\frac{5}{3}}+O(\Delta\sp{\frac{4}{3}}).$$
\end{Theorem}

A slight improvement of this bound gave Goncalves et al. in the recent paper \cite{GoMoPi14}. Moreover, they provide a simple and short proof  
and the upper bound given by Theorem~\ref{Gonc} is almost best possible. \\

\begin{Theorem} \label{Gonc}
(Goncalves, Montassier, Pinlou, 2014, \cite{GoMoPi14}) \\
Let $G$ be a graph with maximum degree $\Delta\geq 3$. Then 
$$\pi(G)\leq \pi_{\l}(G)\leq \left\lceil \Delta\sp{2}+\frac{3}{2\sp{\frac{2}{3}}}\Delta\sp{\frac{5}{3}}+\frac{2\sp{\frac{2}{3}}\Delta\sp{\frac{5}{3}}}{\Delta\sp{\frac{1}{3}}-2\sp{\frac{1}{3}}}\right\rceil=\Delta\sp{2}+\frac{3}{2\sp{\frac{2}{3}}}\Delta\sp{\frac{5}{3}}+O(\Delta\sp{\frac{4}{3}});$$
$\left(\frac{3}{2\sp{\frac{2}{3}}}\dot{=}1.89\right)$. \\
\end{Theorem}

To the set of general results we can also assign a result of Czerwi\' nski and Grytczuk \cite{CzGr07} who proved that for every graph $G$ with maximum degree $\Delta$ there exists a vertex colouring from lists of size at least $16\Delta(G)\sp{2-\frac{1}{k}}$ with no repetitive path on at most $2k$ vertices. 

In \cite{Gr07} various questions concerning nonrepetitive colourings of graphs have been formulated. The open questions from \cite{Gr07} related to $\pi(G)$ and $\pi_{\l}(G)$ will be mentioned throughout this paper too. 

\section{Special classes of graphs} 

There are some classes of graphs, where the Thue chromatic number is known exactly or there are given better upper or lower bounds than the general ones for the graphs belonging to these families. In this section we give an overview of the results on graphs with bounded path-width and tree-width, on planar, outerplanar, cubic, series-parallel, bipartite and complete multipartite graphs.  

\subsection{Planar graphs}
\label{sec:6}

\subsubsection{Paths} 

Thue \cite{Th06} has shown that there are arbitrarily long nonrepetitive sequences over three symbols. As a consequence of this theorem we immediately have a result on nonrepetitive vertex colourings of paths: \\

\begin{Theorem}
(Thue, 1906, \cite{Th06}) \\ 
Let $P_n$ be a path on $n$ vertices. Then $\pi(P_1)=1$, $\pi(P_2)=\pi(P_3)=2$ and for $n>3$ $\pi(P_n)=3$. \\
\end{Theorem}

In the paper \cite{CzGr07} Czerwi\' nski and Grytczuk conjectured that analogue to the Thue theorem  also holds for nonrepetitive list colouring of paths. 

Using the Lov\'asz local lemma Grytczuk, Przyby\l o and Zhu \cite{GPZ10} proved that the assignation of lists of length 4 is always satisfactory for creating nonrepetitive vertex colouring of arbitrary long path using the colours only from the lists preassigned to the vertices of the path: \\
 
\begin{Theorem} \label{GPZ4}
(Grytczuk, Przyby\l o, Zhu, 2011, \cite{GPZ10}) \\
Every path $P_n$ satisfies $\pi_{\l}(P_n)\leq 4$. \\
\end{Theorem}

A more constructive proof of Theorem~\ref{GPZ4} can be found in the paper of Grytczuk, Kozik and Micek, \cite{GKM10}.  \\
Even if the proved bound differs from the conjectured one only by 1, the following problem of Czervi\' nski and Grytczuk remains open: \\

\begin{Question} 
(Czervi\' nski, Grytczuk, 2007, \cite{CzGr07}; Grytczuk, 2007, \cite{Gr07a}) \\
Does every path $P_n$ have a nonrepetitive colouring from arbitrary lists of size three? \\
\end{Question}

The above question is the most interesting open problem from this area, while it is already known that in general the Thue chromatic number and the Thue choice number of the same graph may have arbitrary large difference.

\subsubsection{Cycles} 

A concrete problem on nonrepetitive colourings of cycles was formulated in \cite{AGHR02}. Althrough it concerned the edge variant of the problem, in the family of cycles the same problem can be formulated for vertex colourings. More concretelly, whether $\pi(C_n)=3$ for all $n\geq 18$. 
At the time when the question was asked, it was known that for every cycle of length $n\in\{5$, $7$, $9$, $10$, $14$, $17\}$ $\pi(C_n)=4$, for no other value of $n$ up to 2001 is $\pi(C_n)=4$ and the number of nonrepetitive sequences groves exponencialy with $n$.  The positive answer to this question was given by Currie \cite{Cur02}: \\

\begin{Theorem} 
(Currie, 2002, \cite{Cur02}) \\
For every cycle of length $n\in\{5$, $7$, $9$, $10$, $14$, $17\}$ $\pi(C_n)=4$  and  for other lengths of cycles on at least $3$ vertices $\pi(C_n)=3$. \\
\end{Theorem}

As a corollary of this result we have that every cycle has a subdivision $H$ with $\pi(H) = 3$ (see \cite{Cur02}).

The upper bound for $\pi_{\l}(C_n)$ can be derived from Theorem~\ref{GPZ4}:    \\

\begin{Corollary} 
(Dujmovi\'c, Joret, Kozik, Wood, 2015+, \cite{DJKW13})\\
Every cycle is nonrepetitively $5$-choosable. \\
\end{Corollary} 

To see that consider a cycle $C_n$ with preasign lists of colours of the length at least $5$. Precolour one vertex, remove this colour from every other list and apply the nonrepetitive $4$-choosability result for paths from  \cite{GKM10} or \cite{GPZ10}. \\

Therefore, the following questions are still interesting: \\

\begin{Question} 
(Dujmovi\'c, Joret, Kozik, Wood, 2015+, \cite{DJKW13})\\
Is every cycle nonrepetitively $4$-choosable? \\
\end{Question}

\begin{Question} 
(Dujmovi\'c, Joret, Kozik, Wood, 2015+, \cite{DJKW13})\\
Which cycles are nonrepetitively $3$-choosable? \\
\end{Question}

\subsubsection{Trees} \label{sstrees} 

The first result on the Thue chromatic number for trees was formulated in \cite{AGHR02}, namely, that for every tree $T$ with $\Delta(T) \geq 2$  is $\pi(T) \leq 4$, although it was not proved here. 
The correct proof was given only a few years later by Bre\v sar et al. in \cite{BGKNP07}. \\

\begin{Theorem} \label{tree} 
(Bre\v sar, Grytczuk, Klav\v zar,  Niwczyk, Peterin, 2007, \cite{BGKNP07}) \\ 
If $T$ is a tree, then $\pi(T) \leq 4$, and the bound is tight. \\
\end{Theorem}

Bre\v sar et al. \cite{BGKNP07} also showed a result on Thue chromatic number of subdivisions of trees and a result on Thue chromatic number of trees with small radius. Recall that the {\em eccentricity of a vertex} $u$ is the maximum distance between $u$ and any other vertex, and that the {\em radius of a graph} $G$, denoted $rad(G)$, is the minimum eccentricity of its vertices. \\

\begin{Lemma} 
(Bre\v sar, Grytczuk, Klav\v zar,  Niwczyk, Peterin, 2007, \cite{BGKNP07}) \\ 
Let $T$ be a tree of $rad(T)\leq 4$. 
Then $\pi(T)\leq 3$. \\
\end{Lemma}

\begin{Theorem} 
(Bre\v sar, Grytczuk, Klav\v zar,  Niwczyk, Peterin, 2007, \cite{BGKNP07}) \\ 
Every tree has a subdivision $H$ such that $\pi(H) = 3$. \\
\end{Theorem} 

A family of 4-critical trees is a subfamily of trees for which $\pi(T)=4$. The following question on 4-critical trees is still open:  \\

\begin{Question} 
(Bre\v sar, Grytczuk, Klav\v zar,  Niwczyk, Peterin, 2007, \cite{BGKNP07}) \\ 
Are there infinitely many 4-critical trees? \\
\end{Question} 

Fiorenzi et al. \cite{FOOMZ11} proved that no such result as Theorem~\ref{tree} is possible for nonrepetitive choosability: \\

\begin{Theorem} \label{unbounded} 
(Fiorenzi, Ochem, Ossona de Mendez, Zhu, 2011, \cite{FOOMZ11}) \\ 
For every constant $c$ there is a tree $T$ such that $\pi_{\l}(T)>c$.  \\ 
\end{Theorem}

By this result they gave a negative answer for question of Grytczuk et al. \cite{GKM10} whether the Thue choice number of trees is bounded by a constant. In the same paper \cite{FOOMZ11} Fiorenzi et al. showed the assymptotical behaviour of the Thue choice number of trees of order $n$ and that graphs of bounded tree-depth have bounded Thue choice number. \\
The {\em tree-depth of a graph} can be defined as follows (\cite{NOM06}, \cite{FOOMZ11}): The {\em closure of a rooted tree} $(T , r)$, is defined as the graph $clos(T , r)$ in which $V(clos(T , r)) = V(T )$ and $v_1v_2 \in E(clos(T , r))$ if and only if $v_1$ is an ancestor of $v_2$ or $v_2$ is an ancestor of $v_1$ in $(T , r)$. For a connected graph $G$, the tree-depth of $G$ is the least integer $h$ such that there is a rooted tree $(T , r)$ of height $h$ such that $G$ is a subgraph of $clos(T , r)$. For a disconnected graph $G$, its tree-depth is the maximum of the tree-depth of its connected components. \\

\begin{Theorem}
(Fiorenzi, P. Ochem, P. Ossona de Mendez, X. Zhu, 2011, \cite{FOOMZ11}) \\ 
For every positive integer $h$, the maximum Thue choice number of graphs of tree-depth $h$ is equal to $h$. \\
\end{Theorem}

\begin{Theorem}
(Fiorenzi, P. Ochem, P. Ossona de Mendez, X. Zhu, 2011, \cite{FOOMZ11}) \\ 
The maximum Thue choice number of trees of order $n$ asymptotically satisfies
$\max_{|T|=n} \pi_{\l}(T)=\Omega\left(\left(\frac{\log n}{\log(\log n)}\right)\sp{\frac{1}{2}}\right)$. \\
\end{Theorem}

Kozik and Micek \cite{KM13} proposed an almost linear bound for $\pi_{\l}(T)$ in $\Delta$ for trees: \\

\begin{Theorem}
(Kozik, Micek, 2013, \cite{KM13}) \\
For arbitrary tree $T$ with maximum degree $\Delta$ and for every $\varepsilon > 0$ there is a constant $c$ such that $\pi_{\l}(T)\leq c\cdot\Delta^{1+\varepsilon}$. \\
\end{Theorem} 

Fiorenzi et al. \cite{FOOMZ11} proved that for any $\Delta$ there is a tree $T$ such that $\pi_{\l}(T)=O\left(\frac{\log\Delta}{\log(\log{\Delta})}\right)$. 

In a special case, when the tree a star $S_n$ on $n+1$ vertices, it can be easy observed that $\pi(S_n)=\pi_{\l}(S_n)=2$ (see \cite{PSST14}). \\

Some bounds on Thue chromatic number of trees (e. g. {\em caterpillars} -  trees  in which all the vertices are within distance 1 of a central path) are also in Subsection~\ref{width}.

\subsubsection{Planar graphs in general} 

In \cite{AGHR02} Alon et al. asked a question whether the Thue chromatic number of planar graphs is bounded from above. An equivalent question was posted by Grytczuk: \\

\begin{Question} \label{planarN}
(Alon, Grytczuk, Ha\l uszczak, Riordan, 2002, \cite{AGHR02}; Grytczuk, 2007, \cite{Gr07a}, \cite{Gr07b}, \cite{Gr07}) \\ 
Is there a positive integer $n$ such that $\pi(G)\leq n$ for every planar graph $G$? \\
\end{Question}

Towards Question~\ref{planarN} Bar\'at and Varj\'u \cite{BaVa07} conjectured an upper bound $N=10^{10}$.   \\
Although a similar problem to Question~\ref{planarN} was solved in \cite{BaCz11} for a weaker parameter of nonrepetitive colourings - the facial Thue chromatic number\footnote{Instead of considering that every path in graph is coloured nonrepetitively, only nonrepetitive colouring of every facial path is required.} (see also \cite{HaJe11}, \cite{HaJe12}), Question~\ref{planarN} remains open for the Thue chromatic number of planar graphs under no other condition. 

Dujmovi\'c et al. \cite{DFJW12} asked a question about existence of a logarithmic upper bound for the Thue chromatic number of planar graph with respect to its order.  
In \cite{DFJW12} a logarithmic upper bound of the following form was proved:  

\begin{Theorem} 
(Dujmovi\'c, Frati, Joret, Wood, 2013, \cite{DFJW12}) \\ 
For every planar graph $G$ with $n$ vertices, $\pi(G)\leq 8(1 + \log_{\frac{3}{2}} n)$. \\
\end{Theorem}

Some of the results on Thue chromatic number of planar graphs can be derived from Theorem~\ref{pomocna} too.  

\begin{Theorem} \label{pomocna}
(Dujmovi\'c, Frati, Joret, Wood, 2013, \cite{DFJW12}) \\ 
There is a constant $c$ such that, for every integer $k\geq 1$, every planar graph $G$ is $c\sp{k\sp{2}}$-colourable such that $G$ contains no repetitively coloured path of order at most $2k$. \\
\end{Theorem}

For $k = 2$ Theorem~\ref{pomocna} corresponds to {\em star colourings} - see Subsection~\ref{Not}. These were investigating in \cite{ACKKR04} where Theorem~\ref{chordal} was proved\footnote{In this paper it was also proved that every planar graph is star colourable with $20$ colours.}. \\ 
A graph in which all cycles of four or more vertices have a chord is called {\em chordal graph}. 
The {\em clique number} of a graph $G$, denoted by $\omega(G)$, is the order of a largest complete subgraph of $G$. \\

\begin{Theorem} \label{chordal} 
(Albertson, Chappell, Kierstead, K\"undgen, Ramamurthi, 2004, \cite{ACKKR04}) \\
There is a sequence of chordal graphs $G_1$,$G_2$,$G_3$, . . . such that $\omega(G_t)=t$ 
and $\chi_{st}(G_t)=\frac{(t+1)!}{2(t-1)!}$
Moreover, $G_3$ is outerplanar and $G_4$ is planar. \\
\end{Theorem}

As a corollary of Theorem~\ref{chordal} we have: \\

\begin{Corollary}
(Albertson, Chappell, Kierstead, K\"undgen, Ramamurthi, 2004, \cite{ACKKR04}) \\
There exists a planar graph with $\pi(G)\geq 10$. \\ 
\end{Corollary}

A construction of a planar graph with $\pi(G)= 10$ was found by Bar\'at and Varj\'u \cite{BaVa07}. Ochem showed how to adapt a construction of Albertson et al. \cite{ACKKR04} to the results of Bar\'at and Varj\'u \cite{BaVa07} in order to give an example of a graph $G$ with $\pi(G)= 11$. His result was published only in Appendix of the paper of Dujmovi\'c et al. \cite{DFJW12}: \\ 

\begin{Theorem} 
(Dujmovi\'c, Frati, Joret, Wood, 2013, \cite{DFJW12}) \\ 
There exists a planar graph $G$ with $\pi(G) \geq 11$. \\
\end{Theorem}

Dujmovi\'c et al. \cite{DFJW12} mentioned a class of graphs they consider to be problematic for nonrepetitive colouring and posted up to now still open question regarding the graphs from the family: \\ 
Let $T$ be a tree rooted at a vertex $r$. Let $V_i$ be the set of vertices in $T$ at distance $i$ from $r$. Draw $T$ in the plane with no crossings. Add a cycle on each $V_i$ in the cyclic order defined by the drawing to create a planar graph $G_T$. \\

\begin{Question} 
(Dujmovi\'c, Frati, Joret, Wood, 2013, \cite{DFJW12}) \\
Is $\pi(G_T)\leq c$ for some constant $c$ independent of $T$?  
\end{Question}

\subsection{Outerplanar graphs}
\label{sec:5}
 
At the Budapest workshop in honor of Mikl\'os Simonovits' 60th birthday (June $16^{th}-27^{th}$ 2003), Grytczuk suggested a small change in Question~\ref{planarN}, namely,  replacing planar by outerplanar. 
Independently, Bar\'at and Varj\'u \cite{BaVa07} and K\"undgen and Pelsmajer \cite{KuPe08} found a positive answer to that question: \\

\begin{Theorem}  
(Bar\'at, Varj\'u, 2007, \cite{BaVa07}; K\"undgen and Pelsmajer, 2008, \cite{KuPe08}) \\ 
If $G$ is an outerplanar graph, then $\pi(G) \leq 12$. \\
\end{Theorem} 

A lower bound for the Thue chromatic number of  outerplanar graphs can be derived from Theorem~\ref{chordal}: \\

\begin{Corollary} \label{outerp}
(Albertson, Chappell, Kierstead, K\"undgen, Ramamurthi, 2004, \cite{ACKKR04}) \\
There exists an outerplanar graph with $\pi(G)\geq 6$. \\
\end{Corollary}

A construction of a graph mentioned in Corollary~\ref{outerp} was given by  Bar\'at and Varj\'u \cite{BaVa07}. Moreover, their result is stronger: \\

\begin{Theorem}
(Bar\'at, Varj\'u, 2007, \cite{BaVa07}) \\
There exists an outerplanar graph $G$ with  $\pi(G) \geq 7$. 
\end{Theorem}

\subsection{Graphs with bounded path-width and tree-width} \label{width}

The {\em tree-width}  of a graph G can be defined as the minimum integer $k$ such that $G$ is a subgraph of a chordal graph with no clique on $k + 2$ vertices. Hence, the tree-width of a graph $G$ can be expressed as $\min\{\omega(H)-1;E(G)\subseteq E(H);H$ chordal$\}$, where $\omega(H)$ is the clique number of a graph $H$. 
A {\em path-decomposition} of a graph $G$ is a sequence of subsets of vertices of $G$ such that the endpoints of each edge appear in one of the subsets and such that each vertex appears in a contiguous subsequence of the subsets \cite{RoSe83}.  
The path-width of $G$ is the minimum width of a path decomposition of $G$. 

In \cite{DJKW13} Dujmovi\'c et al. upperbounded the Thue chromatic number of graphs with given path-width $\theta$: \\

\begin{Theorem} 
(Dujmovi\'c, Joret, Kozik, Wood, 2015+, \cite{DJKW13})\\
For every graph $G$ with path-width $\theta$,
$\pi(G)\leq 2\theta\sp{2} + 6\theta + 1$. \\
\end{Theorem}

They supposed that this bound is far from being tight and formulated an open problem whether $\pi(G)\in O(\theta)$ for every graph $G$ with path-width $\theta$. \\
The other question they formulated ask for a relationship between nonrepetitive choosability and path-width. They showed that the graphs with path-width $1$ (i.e., caterpillars) are nonrepetitively $c$-choosable for some constant $c$: \\
 
\begin{Theorem} 
(Dujmovi\'c, Joret, Kozik, Wood, 2015+, \cite{DJKW13}) \\ 
Every caterpillar is nonrepetitively $148$-choosable. \\
\end{Theorem}

A natural question is how the situation looks like for graphs with path-width at least two: \\

\begin{Question} 
(Dujmovi\'c, Joret, Kozik, Wood, 2015+, \cite{DJKW13})\\
Is every graph (or tree) with path-width $2$ nonrepetitively $c$-choosable for some constant $c$? \\
\end{Question}

A class of graphs with bounded tree-width was also investigated.  

We say that a tournament $\textbf{T}$ has the property $S_k$ if and only if any $k$ vertices $v_1;v_2;\dots;v_k$ have a common out-neighbour. Let $f(k)$
be the smallest positive integer such that there exists a tournament on $f(k)$ vertices with property $S_k$. \\
Bar\'at and Varj\'u \cite{BaVa07} proved the following: \\

\begin{Theorem} 
(Bar\'at, Varj\'u, 2007, \cite{BaVa07}) \\ 
Let G be a graph with tree-width at most $k$. Then
$\pi (G) \leq 3\sp{k}\cdot f(k)$. \\
\end{Theorem}

Independently from Bar\'at and Varj\'u \cite{BaVa07}, K\"undgen and Pelsmajer \cite{KuPe08} proved an upper bound for the Thue chromatic number exponential in the tree-width, but independent of the number of vertices: \\

\begin{Theorem}
(K\"undgen and Pelsmajer, 2008, \cite{KuPe08}) \\
If $G$ is a graph of tree-width $k$, then $\pi(G) \leq 4\sp{k}$. \\
\end{Theorem} 

K\"undgen and Pelsmajer also asked a question whether there is a polynomial bound on $\pi(G)$ for graphs of tree-width $k$ \cite{KuPe08}. This was answered in \cite{BaWo05} 
under the additional assumption of bounded degree.  
In particular, Bar\'at and Wood proved an $O(k\Delta)$ upper bound on Thue chromatic number of graph with tree-width $k$: \\ 

\begin{Theorem} 
(Bar\'at, Wood, 2008, \cite{BaWo05})  \\
Every graph $G$ with tree-width $k$ and maximum degree 
$\Delta\geq 1$ satisfies $\pi(G)\leq 10(k+1)(\frac{7}{2}\Delta-1)$. \\
\end{Theorem}

Another question was asked by Dujmovi\'c et al. \cite{DJKW13}: \\

\begin{Question} 
(Dujmovi\'c, Joret, Kozik, Wood, 2015+, \cite{DJKW13}) \\
Is $\pi(G)$ bounded from above by a polynomial function of tree-width of the graph $G$? \\
\end{Question} 

This is still open. Recall, that the tree-width of a graph does not provide an upper bound on its Thue choice number $\pi_{\l}(G)$ - see Subsection~\ref{sstrees}. \\

The lower bound for the Thue chromatic number of graphs with given tree-width $k$ comes from the theorem of Albertson et al. \cite{ACKKR04}: \\

\begin{Theorem} 
(Albertson, Chappell, Kierstead, K\"undgen, Ramamurthi, 2004, \cite{ACKKR04})\\
There exists a graph $G$ with tree-width $k$ and $\pi(G)\geq \chi_{st}(G)=\frac{(k+2)!}{2\cdot k !}$. 
\end{Theorem}

\subsection{Other special classes of graphs}
 
A {\em cubic graph} is a graph in which all vertices have degree three. \\
Grytczuk \cite{Gr06} asked a question how large the Thue chromatic number for cubic graphs can be.  
Using the probabilistic approach he proved the following: \\

\begin{Theorem}
(Grytczuk, 2006, \cite{Gr06}) \\ 
Let $G$ be a cubic graph. Then $\pi(G)\leq 108$. \\
\end{Theorem}

A graph is a {\em series-parallel graph} (see \cite{Duff65}), if it may be turned into $K_2$; $V(K_2)=\{s, t\}$, $E(K_2)=e$, by a sequence of the following operations: \\
1. Replacement of a pair of parallel edges with a single edge that connects their common endpoints. \\
2. Replacement of a pair of edges incident to a vertex of degree 2 other than s or t with a single edge. 
 
Bar\'at and Varj\'u showed the following bound on $\pi(G)$ for a series-parallel graph:  \\
 
\begin{Theorem} 
(Bar\'at, Varj\'u, 2007, \cite{BaVa07}) \\ 
Let $G$ be a series-parallel graph. Then $\pi(G) \leq 63$. \\
\end{Theorem}

A graph that does not contain any odd-length cycles is called {\em bipartite graph}. One result on Thue chromatic number of bipartite graphs comes from a construction of K\"undgen and Pelsmajer \cite{KuPe08}: \\

\begin{Theorem}
(K\"undgen and Pelsmajer, 2008, \cite{KuPe08}) \\
There are bipartite graphs with arbitrarily high girth and Thue chromatic number. \\
\end{Theorem}  

By this result they gave a negative answer for the question of Schaefer and Umans whether the Thue chromatic number of arbitrary graph $G$ can be bounded from above by some absolute constant $k$ \cite{SchUm02}. \\

An {\em independent set} of vertices in a graph is a set of vertices where no two vertices are adjacent. 
The {\em independence number} of $G$, $\alpha(G)$, is the size of the largest independent set of a given graph G.
For an integer $k\geq 2$ and for positive integers $n_1$, $n_2,\dots$, $n_k$ a {\em complete $k$-partite graph} $K_{n_1,n_2,\dots,n_k}=G(V,E)$ is a graph whose vertex set $V(G)$ can be partitioned into $k$ independent sets $V_1$, $V_2,\dots,$ $V_k$, with $|V_i|=n_i$ for $i=1,2,\dots,k$, such that $u,v\in E(G)$ if $u\in V_i$ and $v\in V_j$, where $1\leq i$, $j\leq k$ and $i\neq j$. A {\em complete multipartite graph} is a graph that is complete $k$-partite for some $k$ \cite{ChZh09}.

Peterin et al. \cite{PSST14} proved the following: \\

\begin{Theorem} \label{nasa}
(Peterin, Schreyer, \v Skrabu\ml \'akov\'a,  Taranenko, 2014, \cite{PSST14}) \\
Let $G$ be a graph on $n$ vertices and $\alpha(G)$ be an independence number of $G$. Then $\pi (G)\leq \pi_{\l}(G)\leq n-\alpha (G)+1$. \\
If $G$ is a complete multipartite graph, then $\pi (G)=\pi_{\l}(G)=n-\alpha (G)+1$. \\
\end{Theorem}

A corollary of Theorem~\ref{nasa} is that in a case of complete bipartite graphs, as well as complete graphs, the Thue chromatic number equals the Thue choice number: \\  

\begin{Corollary}
(Peterin, Schreyer, \v Skrabu\ml \'akov\'a,  Taranenko, 2014, \cite{PSST14}) \\
$\pi(K_n)=\pi_{\l}(K_n)=n$ for the complete graph $K_n$ on $n$ vertices. \\
$\pi(K_{m,n})=\pi_{\l}(K_{m,n})=\min\{m,n\}+1$ for a complete bipartite graph $K_{m,n}$. 
\end{Corollary}

\section{Products of graphs}
\label{sec:9}

A square grid can be understand as a Cartesian product of two path graphs. \\
In general the {\em Cartesian product} $G \square H$ of graphs $G$ and $H$ is a graph such that the vertex set of $G \square H$ is the Cartesian product $V(G)\times V(H)$; and two vertices $(v_1,w_1)$ and $(v_2,w_2)$ are adjacent in $G \square H$ if and only if either $v_1 = v_2$ and $w_1$ is adjacent with $w_2$ in graph $H$, or $w_1 = w_2$ and $v_1$ is adjacent with $v_2$ in $G$. \\

Bar\'at and Varj\'u \cite{BaVa07} and independently K\"undgen and Pelsmajer \cite{KuPe08} proved that the $k \times k$ square grid has a bounded Thue chromatic number: \\

\begin{Theorem} 
(Bar\'at, Varj\'u, 2007, \cite{BaVa07}; K\"undgen, Pelsmajer, 2008, \cite{KuPe08}) \\
Every square grid graph admits a nonrepetitive $16$-colouring. \\
\end{Theorem}

The {\em strong product} of graphs $G$ and $H$, $G\boxtimes H$, is the graph with vertex set $V(G)\times V(H)$ in which distinct vertices $(v_1,w_1)$ and $(v_2,w_2)$ are adjacent when $\deg_G(v_1,v_2)\leq 1$ and $\deg_H(w_1, w_2)\leq 1$. \\
The upper bound for the Thue chromatic number of the strong product of two path graphs was showed by K\"undgen and Pelsmajer \cite{KuPe08}: \\

\begin{Theorem}
(K\"undgen, Pelsmajer, 2008, \cite{KuPe08}) \\
The strong product of $t$ paths admits a nonrepetitive colouring with at most $4\sp{t}$ colours. \\
\end{Theorem}

The {\em lexicographic product} $G[H]$, or {\em blow-up} of $G$ by $H$, of graphs $G=(V_1,E_1)$ and $H=(V_2,E_2)$ is a graph such that the vertex set of $G[H]$ is the Cartesian product $V(G)\times V(H)$ and two vertices $(v_1,w_1)$ and $(v_2,w_2)$ are adjacent in $G[H]$ if either $v_1$ is adjacent with $v_2$ in $G$ 
or $v_1=v_2$ and $w_1$ is adjacent with $w_2$ in $H$.  

The lexicographic product of graphs was first studied by Felix Hausdorf in 1914 - see \cite{Ha14} for the original book or \cite{Ha65} for its reprint.  But nonrepetitive colourings of lexicographic product of graphs have been systematicly studied only recently. However, some of the older results can be also transformed into words of nonrepetitive colouring of lexicographic product of graphs. One of these results was achieved by Bar\'at and  Wood: \\

\begin{Theorem} 
(Bar\'at, Wood, 2008, \cite{BaWo05}) \\
For any tree $T$ and integer $k$, $\pi(T[K_k])\leq 4k$. \\
\end{Theorem}

This bound is tight, as for every positive integer $k$ there exist a tree $T$ for which $\pi(T[E_k])=4k$ (see \cite{KPZ13}).

The Thue chromatic number of $G[H]$ when $G$ is a path and $H$ is either an empty graph $E_k$ on $k$ vertices or a complete graph $K_k$ was studied in \cite{KPZ13}. The main results of Keszegh et al. \cite{KPZ13} are the following: \\

\begin{Theorem} 
(Keszegh, Patk\'os, Zhu, 2013, \cite{KPZ13}) \\ 
For every $n\geq 4$ and $k\neq 2$, $\pi(P_n[E_k])=2k+1$. For $k=2$, $5\leq\pi(P_n[E_2])\leq 6$. \\
\end{Theorem}

\begin{Theorem} 
(Keszegh, Patk\'os, Zhu, 2013, \cite{KPZ13}) \\
For every integer $n\geq 28$, $3k+\lfloor \frac{k}{2}\rfloor\leq \pi(P_n[K_k])\leq 4k$. \\
\end{Theorem}

By further requirement that every copy of $E_k$ is rainbow-coloured Keszegh et al. \cite{KPZ13} proved that the smallest number of colours needed for $P[E_k]$ is at least $3k+1$ and at most $3k+\lceil\frac{k}{2}\rceil$. \\
In the case when $G$ is outerplanar graph they proved that $\pi(G[K_k])\leq 16k$ for any integer $k\geq 2$. Furthermore, for every positive integer $k$ there exists an outerplanar graph $G_0$ such that $\pi(G_0[E_k])>6k$.

In Keszegh et al. \cite{KPZ13} are formulated some open questions on Thue chromatic number of lexicographic product of graphs: \\

\begin{Question} 
(Keszegh, Patk\'os, Zhu, 2013, \cite{KPZ13}) \\
Is there a constant $c$ such that $\pi_{\l}(P_\infty[K_k])\leq c\cdot k$? \\
\end{Question}

\begin{Question} 
(Keszegh, Patk\'os, Zhu, 2013, \cite{KPZ13}) \\
Is there a function $f$ such that for every graph $G$ of maximum degree $\Delta$, $\pi(G[K_k])\leq k\cdot f(\Delta)$? Perhaps $\pi(G[K_k])\leq c\cdot k\Delta\sp{2}$ for some constant $c$? \\
\end{Question}

Peterin et al. \cite{PSST14} asked whether the Thue chromatic number of lexicographic product of graphs $G$ and $H$, $\pi(G[H])$, can be bounded from below by a function of $\pi (H)$ and $\pi(G)$: \\

\begin{Question}
(Peterin, Schreyer, \v Skrabu\ml \'akov\'a,  Taranenko, 2014, \cite{PSST14}) \\
Is it true that for all simple graphs $G$ and $H$ we have that $\pi (H)+(\pi (G)-1)|V(H)|\leq \pi (G[H])$ ? \\
\end{Question}

According to the results of Peterin et al. \cite{PSST14} the conjecture is true in the case when $G$ is a complete multiparite graph and $H$ is an arbitrary graph.  
Moreover, in this case it holds: $\pi(H)+(\pi(G)-1)|V(H)|= \pi(G[H]) = \pi(H)+(|V(G)|-\alpha (G))|V(H)|$,     
because $\pi (G[H])\leq \pi (H)+(|V(G)|-\alpha (G))|V(H)|$ holds for all simple graphs $G$ and $H$ \cite{PSST14}. 

Using these results together with Theorem~\ref{nasa} another exact bounds for the Thue chromatic number of $G[H]$ can be derived: \\

\begin{Corollary}
(Peterin, Schreyer, \v Skrabu\ml \'akov\'a,  Taranenko, 2014, \cite{PSST14}) \\
For any graph $H$ it holds  \\
$\pi(K_n[H])=\pi(H)+(n-1)|V(H)|$, where $K_n$ is complete graph on $n$ vertices. \\
$\pi(S_n[H])=\pi(H)+|V(H)|$, where $S_n$ is a star on $n+1$ vertices. \\
$\pi(K_{m,n}[H])=\pi(H)+\min\{m,n\}\cdot |V(H)|$, where $K_{m,n}$ is a complete bipartite graph on $m+n$ vertices. 
\end{Corollary}

\section{Subdivisions of graphs}
\label{sec:10}

An easy corollary of the Thue theorem \cite{Th06} is that every path has a subdivision $P_S$ with $\pi(P_S) = 3$. 
A natural question then is whether every graph has a subdivision that is nonrepetitively $3$-colourable. This was formulated by Bre\v sar et al. \cite{BGKNP07}. 

Grytczuk \cite{Gr07a} proved that every graph $G$ has a subdivision $G_S$ with $\pi(G_S)\leq 5$.
Bar\'at and Wood \cite{BaWo05} improved this result by showing that every graph $G$ has a subdivision $G_S$ with $\pi(G_S)\leq 4$. Indepedently Marx and Schaefer \cite{MaS09} proved the same upper bound on $\pi(G_S)$. Using the similar approach, finally, Pezarski and Zmarz \cite{PeZm09} confirmed that every graph $G$ has a subdivision $G_S$ such that $\pi(G_S)\leq 3$.   
Moreover, they derived a formula that gives the exact value of $\pi(G_s)$ for a subdivision $G_s$ of any graph $G = (V,E)$: \\

\begin{Theorem} 
(Pezarski, Zmarz, 2009,  \cite{PeZm09}) \\ 
$\pi (G_s)=1$ if $E=\emptyset$,\\
$\pi (G_s)=2$ if $E\neq\emptyset$ and $G$ is a star forest\footnote{any acyclic graph that does not contain $P_4$},\\
otherwise, $\pi (G_s)=3$. \\
\end{Theorem}

Grytczuk \cite{Gr07a} asked also for the bounds on Thue chromatic number of subdivisions of graphs under the restriction for number of vertices subdividing each edge.  \\

\begin{Question} \label{division}
(Grytczuk, 2007, \cite{Gr07a}) \\ 
Are there constants $k$ and $n$ such that every planar graph has a subdivision, with at most $k$ vertices subdividing an edge, which is nonrepetitively $n$-colourable? \\
\end{Question}


Ne\v set\v ril et al. \cite{NOMW11} proved that every graph has a nonrepetitively $17$-colourable subdivision with $O(\log n)$ division vertices per edge, and that $(\log n)$ division vertices are needed on some edge of a nonrepetitively $O(1)$-colourable subdivision of $K_n$: \\

\begin{Theorem} 
(Ne\v set\v ril, Ossona de Mendez, Wood, 2012, \cite{NOMW11}) \\
The $\lceil \log n\rceil$-subdivision of $K_n$ has a nonrepetitive $17$-colouring. \\
Moreover, if $K_{nS}$ is a subdivision of $K_n$ and $\pi(K_{nS})\leq c$, then some edge of $K_n$ is subdivided at least $\log_{c+3} \left(\frac{n}{2}\right) -1$ times. \\
\end{Theorem}

Via results of Ne\v set\v ril et al. \cite{NOMW11} $\pi(G)$ is strongly topological and it is a function of $\pi(G_S)$ and the number of vertices subdividing edges of the graph $G$: \\ 
 
\begin{Theorem} \label{nesetril}
(Ne\v set\v ril, Ossona de Mendez, Wood, 2012, \cite{NOMW11}) \\
There is a function $f$ such that $\pi(G)\leq f (\pi(G_S), d)$ for every ($\leq d$)-subdivision $G_S$ of a graph $G$. \\
\end{Theorem}

Theorem~\ref{nesetril} implies that to prove that planar graphs have bounded Thue chromatic number, it suffices to show that every planar graph $G$ has a subdivision $G_S$ with bounded $\pi(G_S)$ and a bounded number of division vertices per edge. This shows that the Question~\ref{division} and Question~\ref{planarN} are equivalent.

In \cite{DJKW13} a similar question for a general graph $G$ can be found: \\

\begin{Question}
(Dujmovi\'c, Joret, Kozik, Wood, 2015+, \cite{DJKW13}) \\ 
Is there a function $f$ such that every graph $G$ has a nonrepetitively $O(1)$-colourable subdivision with $f(\pi(G))$ division vertices per edge? \\
\end{Question}

In \cite{NOMW11} one can find a lot of results concerning nonrepetitive colourings of subdivided graphs. Among all we mention one more: \\

\begin{Theorem} 
(Ne\v set\v ril, Ossona de Mendez, Wood, 2012, \cite{NOMW11}) \\
For $d\geq 2$, and $K_{nS_d}$ being a ($\leq d$)-subdivision of $K_n$ it holds $\left(\frac{n}{2}\right)^\frac{1}{d+1} \leq \pi(K_{nS_d}) \leq 9 \lceil n^{\frac{1}{n+1}}  \rceil.$  \\ 
\end{Theorem}

Nonrepetitive choosability of subdivided graphs was studied in \cite{DJKW13}, where it was proved that every graph has a nonrepetitively $5$-choosable subdivision. By this result Dujmovi\'c et al. \cite{DJKW13} gave a positive answer to the question of Grytczuk et al. \cite{GPZ10} whether there exist a constant $c$ such that every graph $G$ has a subdivision $G_S$ such that $\pi_l(G_S)\leq c$: \\

\begin{Theorem} \label{ch5}
(Dujmovi\'c, Joret, Kozik, Wood, 2015+, \cite{DJKW13}) \\ 
Let $G_S$ be a subdivision of a graph $G$, such that each edge $vw\in E(G)$ is subdivided at least $\lceil 10\sp{5} \log_2(deg(v)+1)\rceil+\lceil 10\sp{5} \log_2(deg(w)+1)\rceil +2$ times. Then $\pi_{\l}(G_S)\leq 5$. \\
\end{Theorem}

Theorem~\ref{ch5} was proved by the entropy compression method. A similar theorem with more colours and $O(\log\Delta(G))$ division vertices per edges can be proved using the Lov\'asz local lemma: \\

\begin{Theorem} 
(Dujmovi\'c, Joret, Kozik, Wood, 2015+, \cite{DJKW13}) \\ 
For every graph $G$ with maximum degree $\Delta$, every subdivision $G_S$ of $G$ with at least $3+400\log\Delta$ division vertices per edge is nonrepetitively $23$-choosable. \\
\end{Theorem}

Dujmovi\'c et al. \cite{DJKW13} supposed that the upper bound for the Thue choice number for subdivisions of graphs given by Theorem~\ref{ch5} is not best possible. Therefore, they asked the following question: \\

\begin{Question}
(Dujmovi\'c, Joret, Kozik, Wood, 2015+, \cite{DJKW13}) \\ 
Does every graph have a nonrepetitively $4$-choosable subdivision? Even $3$-choosable might be possible.  
\end{Question}

\section{Remarks on the compexity of nonrepetitive colourings}
\label{sec:11}

Marx and  Schaefer \cite{MaS09} showed that deciding whether a colouring is repetitive is $NP$-complete: \\

\begin{Theorem} 
(Marx, Schaefer, 2009, \cite{MaS09}) \\  
Determining whether a particular colouring of a graph is nonrepetitive is coNP-hard, even if the number of colours is limited to four. \\
\end{Theorem}

Marx and  Schaefer \cite{MaS09} also gave an algorithm that is able to check whether $G$ has a repetitive sequence of length $2k$: \\

\begin{Theorem}
(Marx, Schaefer, 2009, \cite{MaS09}) \\  
Given a vertex-coloured graph $G=G(V,E)$, it can be checked in time $k\sp{O(k)}\cdot |V|\sp{5} \log |V|$ whether $G$ has a repetitive sequence of length $2k$. \\
\end{Theorem}

For $k=2$ we get a star-free colouring of graphs without repetitive sequences of length at most $4$. Deciding whether a graph has a star-free colouring with three colours is $NP$-complete, even if the graph is bipartite  \cite{CoMo84}.

A {\em cograph} is any $P_4$-free graph. Some superclasses of cographs are e.g. {\em $P_4$-tidy graphs} and {\em $(q,q-4)$-graphs}. 
A graph is called $(q, q-4)$-graph if no set of at most $q$ vertices induces more than $q-4$ distinct $P_4$ 's \cite{BaOl98}. 
A graph is called $P_4$-tidy if for every $P_4$ induced by $(u, v, x, y)$, there exists at most one vertex $z$ such that $u, v, x, y, z$ induces more than one $P_4$ \cite{GiRoTh97}. 

Lyons \cite{Ly11} obtained a polynomial time algorithm to find an optimal acyclic colouring\footnote{{\em Acyclic colouring} is a proper colouring such that every pair of colour classes induces a forest (see e.g. \cite{SeVo13}).}  and an optimal star colouring of a cograph. In \cite{CLSMMS11} and \cite{CDS12} it was proved that every acyclic colouring of a cograph is also nonrepetitive. Moreover, Campos et al. \cite{CLSMMS11} and Costa et al. \cite{CDS12} showed that there exist linear time algorithms to obtain $\pi(G)$ for $G$ being a $P_4$-tidy or a $(q, q-4)$-graph for some fixed integer $q$. 

Dujmovi\'c et al. \cite{DJKW13} asked also an interesting question regarding graph algorithms: \\

\begin{Question}
(Dujmovi\'c, Joret, Kozik, Wood, 2015+, \cite{DJKW13}) \\
Is there a polynomial-time Monte Carlo algorithm that nonrepetitively $O(\Delta\sp{2})$-colours a graph with maximum degree $\Delta$? \\
\end{Question}

The closest related result was proved by Haeupler et al. \cite{HSS11}: \\

\begin{Theorem}
(Haeupler, Saha, Srinivasan, 2011, \cite{HSS11}) \\
There exists a constant $c>0$ such that for every constant 
$\varepsilon > 0$ there exists a Monte Carlo algorithm that given a graph $H$ with maximum degree $\Delta$, produces a nonrepetitive colouring using at most $c \Delta^{2+\varepsilon}$ colours. The failure probability of the algorithm is an arbitrarily small inverse polynomial in the size of $H$. 
\end{Theorem} 


\section{Conclusion}
\label{sec:13}

In this section we give a summary of some results presented in previous sections by creating comparating tables on the values of Thue chromatic number and Thue choice number of selected families of graphs. 

The results comparing $\pi(G)$ and $\pi_{\l}(G)$ for selected subfamilies of graphs can be found in Table~1. 

\begin{table}[h!] 
\caption{Comparation of $\pi(G)$ and $\pi_{\l}(G)$ for graph $G$ from given family of graphs}
		\begin{center}
		\begin{small}
		\begin{tabular}[t]{|l|c|c|} 
		\hline
{\bf Graph $G$}& \ {\bf Result on $\pi(G)$} & \ {\bf Result on $\pi_{\l}(G)$} \\ 
\hline
path $P_n$ on $n$ vertices & $\pi(P_n)=3$ for $n>3$ & $\pi_{\l}(P_n)\leq 4$ \\ \hline
cycle $C_n$ on $n$ vertices & $\pi(C_n)=3$ for $n\notin M$ & $\pi_{\l}(C_n)\leq 5$ \\ 
  & $\pi(C_n)=3$ for $n\in M$ &   \\ 
    & $M=\{5, 7, 9, 10, 14, 17\}$ &   \\\hline 
star $S_n$ on $n+1$ vertices & $\pi(S_n)=2$ & $\pi_{\l}(S_n)=2$ \\ \hline
caterpillar $H$ & $\pi(H)\leq 4$ & $\pi_{\l}(H)\leq 148$ \\ \hline
tree $T$ of maximum degree $\Delta$ & $\forall T: \pi(T)\leq 4$ & $\forall \varepsilon >0 \ \exists c\in \mathbb{R}:\pi_{\l}(T)\leq c\Delta^{1+\varepsilon}$ \\
  & $\exists T:\pi(T)=4$ & $\forall c\in \mathbb{R} \ \exists T:\pi_{\l}(T)\geq c$ \\ \hline
planar graph $G$ & $\exists G:\pi(G)\geq 11$ & $\exists G:\pi_l(G)\geq 11$ \\ \hline	
outerplanar graph $G$ & $\exists G: \pi(G)\geq 7$ & $\exists G: \pi_{\l}(G)\geq 7$ \\
 & $\forall G: \pi(G)\leq 12$ &   \\ \hline
$G$ with path-width $\theta$ & $\pi(G)\leq 2\theta\sp{2} + 6\theta + 1$ &  for $\theta = 1$; $\pi_{\l}(G)\leq c$ \\ \hline
complete graph $K_n$ & $\pi(K_n)=n$ &  $\pi_{\l}(K_n)=n$ \\ \hline
complete multipartite graph $G$ & $\pi(G)=n-\alpha(G)+1$ &  $\pi_{\l}(G)=n-\alpha(G)+1$ \\ \hline
complete bipartite graph $K_{m,n}$ & $\pi(K_{m,n})=$min$\{m,n\}+1$ &  $\pi_{\l}(K_{m,n})=$min$\{m,n\}+1$ \\ \hline
bipartite graph $G$ & $\forall c\in\mathbb{R} \ \exists G: \pi(G)\geq c$ &  $\forall c\in\mathbb{R} \ \exists G: \pi_{\l}(G)\geq c$ \\ \hline 
subdivision $G_S$ of a graph $G$ & $\forall G \ \exists G_S: \pi(G_S)\leq 3$ &  $\forall G \ \exists G_S: \pi_{\l}(G_S)\leq 5$ \\
		\hline   
		\end{tabular}
		\end{small}
		\end{center}
	\end{table}	 

As most of the results in the Table~2 were achieved via probabilistic approach, the upper bounds on $\pi(G)$ and $\pi_{\l}(G)$ are here the same:  

\begin{table}[h!]  
\caption{Comparation of $\pi(G)$ and $\pi_{\l}(G)$ for arbitrary graph $G$} 
		\begin{center}
		\begin{small}
		\begin{tabular}[t]{|l|c|} 
		\hline
{\bf Graph $G$}& \ {\bf Result on $\pi(G)$ and $\pi_{\l}(G)$} \\ 
\hline	
graph $G$ & $\exists G: \pi_{\l}(G)\geq \pi(G) \geq c\frac{\Delta^2}{\log{\Delta}}$  \\ \hline
graph $G$  &  $\forall G: \pi(G)\leq \pi_{\l}(G) \leq 10,4 \Delta^2$ \\ 
 &  $\pi(G)\leq \pi_{\l}(G) \leq 12,92 (\Delta-1)^2$ \\ 
 & $\pi(G)\leq \pi_{\l}(G) \leq \Delta\sp{2}+\frac{3}{2\sp{\frac{2}{3}}}\Delta\sp{\frac{5}{3}}+O(\Delta\sp{\frac{4}{3}})$ \\
		\hline 
		\end{tabular}
		\end{small}
		\end{center}
	\end{table}

\textbf{Acknowledgements:} 

This work was supported by the Slovak Research and Development Agency
under the contract No. APVV-14-0892, this work was also supported by the Slovak Research and Development Agency under the contract No. APVV-0482-11, by the grants VEGA 1/0529/15, VEGA 1/0908/15 and KEGA 040TUKE4/2014.

\end{document}